\documentclass[11pt]{article}
\usepackage[latin1]{inputenc}
\usepackage{amsmath,amsthm,amssymb}

\newtheorem{thm}{Theorem}
\numberwithin{thm}{section}

\def\eq#1{(\ref{#1})}

\newcommand{\neweq}[1]{\begin{equation}\label{#1}}

\def\ep{\varepsilon}

\def\phi{\varphi}
\def\RR{\mathbb R}

\def\di{\displaystyle}
\def\ri{\rightarrow}
\def\intom{\int_\Omega}
\def\huo{H^1_0(\Omega )}
\def\incep{\left\{\begin{array}{cl} }
 \def\termin{\end{array}\right. }
\def\2af{2^*_\alpha}

\def\plpm{$(P_\lambda)_\pm$}
\def\plp{$(P_\lambda)_+$}
\def\plm{$(P_\lambda)_-$}
\def\cale{{\mathcal E}}
\def\calf{{\mathcal F}_\lambda}

\title{\sc Combined effects in nonlinear problems arising in the
study of anisotropic continuous media}
\author{\sc Vicen\c tiu R\u adulescu$\,^{a,b}$ and Du\v{s}an Repov\v{s}$\,^{c,d}$\\
\small $^a\,$Institute of Mathematics ``Simion Stoilow" of the Romanian
Academy, Bucharest\\
\small $^b\,$Department of Mathematics,
University of Craiova,  200585 Craiova, Romania\\
 \small $^c\,$Faculty of Mathematics and
Physics, University of Ljubljana,\\ \small Jadranska  19,  P. O. Box 2964, 1001 Ljubljana, Slovenia\\
\small $^d\,$Faculty of Education, University of Ljubljana, Kardeljeva plo\v{s}\v{c}ad 16, 1000 Ljubljana, Slovenia\\
 \small E-mail: {\tt vicentiu.radulescu@imar.ro}\qquad {\tt dusan.repovs@guest.arnes.si}\\}

\date{}

\begin{document}

\maketitle
\begin{center}{\it This paper is dedicated with esteem to Professor Marius Iosifescu on his 75th birthday}\end{center}

\begin{abstract}
We are concerned with the Lane-Emden-Fowler equation $-\Delta u=\lambda k(x)u^q\pm h(x)u^p$ in $\Omega$, subject to the Dirichlet boundary condition $u=0$ on $\partial\Omega$, where $\Omega$ is a smooth bounded domain in $\RR^N$, $k$ and $h$ are variable potential functions, and $0<q<1<p$. Our analysis combines monotonicity methods with variational arguments.
 \\
{\bf Keywords:} Lane-Emden-Fowler equation, bifurcation problem, anisotropic continuous media, perturbation,
positive solution.\\
{\bf 2010 Mathematics Subject Classification:} 35A05, 35B40, 35J60,
58E05.
\end{abstract}

\section{Introduction and the main results}
In this paper we are concerned with the study of combined effects in a class of nonlinear bifurcation problems that arise in anisotropic continuous media. Bifurcation problems have a long history and their treatment goes
back to the XVIIIth century. To the best of our knowledge, the first bifurcation problems
is related to the buckling of a thin rod under thrust and was
investigated by Daniel Bernoulli and Euler around 1744. In the
case in which the rod is free to rotate at both end points, this
yields the one-dimensional bifurcation problem
$$\left\{\begin{array}{ll}
& u''+\lambda\sin u=0\qquad\mbox{in $(0,L)$}\\
& 0\leq u\leq\pi\\
& u'(0)=u'(L)=0\,.\end{array}\right.$$

Consider the bifurcation problem
\begin{equation}\label{1}
\left\{\begin{array}{lll}
&\di -\Delta u=\lambda f(u) & \qquad \mbox{in}\ \Omega \\
&\di u=0 &\qquad \mbox{on}\ \partial\Omega\, ,\\
\end{array}\right.
\end{equation}
where $\Omega\subset\RR^N$ is a bounded domain with smooth boundary. Problems of this type have been studied starting with the pioneering paper \cite{gel} by Gelfand, who considered the case $f(u)=e^u$. We also refer to the important contributions due to Amann \cite{amsiam} and Keller \& Cohen \cite{kelcoh}, who assumed that $f$ is a function of class $C^1$  which is
  convex and positive, such that $f'(0)>0$. They proved the following basic facts:

  (i) there exists $\lambda^*\in(0,\infty)$ such that problem \eq{1} has (resp., has no) classical solution if $\lambda\in (0,\lambda^*)$ (resp., $\lambda>\lambda^*$);

(ii) for  $\lambda\in (0,\lambda^*)$, among the solutions of \eq{1} there exists a minimal one, say $u(\lambda)$;

(iii) the mapping $\lambda\longmapsto u(\lambda)$ is convex, increasing and of class $C^1$;

(iv) $u(\lambda)$ is the only solution of problem \eq{1} such that the operator $-\Delta -\lambda f'(u)$ is coercive.

Motivated by a problem raised by Brezis, Mironescu and R\u adulescu studied in \cite{mirradcras,mirradna} the case where $f$ fulfills the above hypotheses and has a linear growth at infinity, that is, there exists $\lim_{t\ri\infty}f(t)/t=:a\in (0,\infty)$. The results in this framework are strongly related with the sign of $b:=\lim_{t\ri\infty}(f(t)-at)$. For instance , if $b\geq 0$ then the following properties hold: (i) $\lambda^*=\lambda_1/a$; (ii) $u(\lambda)$ is the unique solution of problem \eq{1} for all $\lambda\in (0,\lambda^*)$; (iii) $\lim_{\lambda\ri\lambda^*-}u(\lambda)=+\infty$ uniformly on compact subsets of $\Omega$; and (iv) problem \eq{1} has no solution if $\lambda=\lambda^*$. A different behaviour holds if $b\in [-\infty,0)$.

The perturbed nonlinear problem
\begin{equation}\label{2}
\left\{\begin{array}{lll}
&\di -\Delta u=\lambda f(u)+a(x)g(u) & \qquad \mbox{in}\ \Omega \\
&\di u>0 & \qquad \mbox{in}\ \Omega \\
&\di u=0 &\qquad \mbox{on}\ \partial\Omega\, ,\\
\end{array}\right.
\end{equation}
was studied in \cite{cgrjmpa},
where $f$ has either a sublinear or a linear growth at infinity, $g$ is a singular nonlinearity, and $a$ is a nonnegative potential. The framework discussed in \cite{cgrjmpa} includes the case $f(t)=t^p$ and $g(t)=t^{-\gamma}$, where $0<p\leq 1$ and $0<\gamma<1$. For instance, if $f$ is sublinear and $\inf_{x\in\Omega}a(x)>0$ then problem \eq{2} has a unique solution for all $\lambda\in\RR$.
From
a physical point of view, problem \eq{2} arises in the context of chemical heterogenous
catalysis, in the theory of heat conduction in electrically conducting materials, as well
as in the study of non-Newtonian fluids, boundary layer phenomena for viscous fluids.
Such equations are also encountered in glacial advance (see \cite{r32}), in
transport of coal slurries down conveyor belts (see \cite{r4}), and in several other geophysical and
industrial contents (see \cite{r3} for the case of the incompressible flow of a uniform stream past
a semi-infinite flat plate at zero incidence).

The above results are summarized and completed in the recent works \cite{bare}, \cite{groxford}, \cite{kmrcup}, and \cite{hindawi}.

In this paper we are concerned with the competition between convex and concave nonlinearities and variable potentials. Such problems arise in the study of anisotropic continuous media.
We point out that nonlinear elliptic equations with convex--concave nonlinearities have been studied starting with the seminal paper by Ambrosetti, Brezis and Cerami \cite{ARC}. They considered the Dirichlet problem
\begin{equation}\label{abc}
\left\{\begin{array}{lll}
&\di -\Delta u=\lambda u^{q}+u^{p}, \qquad &\mbox{in}\ \Omega\\
&\di u>0, \qquad&\mbox{in}\ \Omega\\
&\di u=0, \qquad&\mbox{on}\ \partial\Omega\,,
\end{array}\right.
\end{equation}
where $\lambda$ is a positive parameter, $\Omega\subset\RR^N$ is a bounded domain with smooth boundary, and $0<q<1<p<2^\star -1$ ($2^\star=2N/(N-2)$ if $N\geq 3$, $2^\star=+\infty$ if $N=1,\, 2$). Ambrosetti, Brezis and Cerami proved that there exists $\lambda^*>0$ such that problem \eq{abc} admits at least two solutions for all $\lambda\in (0,\lambda^*)$, has one solution for $\lambda=\lambda^*$, and no solution exists provided that $\lambda>\lambda^*$.
In \cite{AT2}, Alama and Tarantello studied the related Dirichlet problem with indefinite weights
\begin{equation}\label{alama}
\left\{\begin{array}{lll}
-\Delta u-\lambda u=k(x)u^q-h(x)u^p, &\mbox{if}& x\in\Omega\\
u>0, &\mbox{if}& x\in\Omega\\
u=0, &\mbox{if}& x\in\partial\Omega\,,
\end{array}\right.
\end{equation}
where $h$, $k$ are nonnegative and $1<p<q$. We refer to \cite{lincei10} for a related problem with lack of compactness.

In the present paper we are concerned with the nonlinear elliptic problem
$$
\left\{\begin{array}{lll}
-\Delta u=\lambda k(x)u^q\pm h(x)u^p, &\mbox{if}& x\in\Omega\\
u>0, &\mbox{if}& x\in\Omega\\
u=0, &\mbox{if}& x\in\partial\Omega\,,
\end{array}\right.\eqno(P_\lambda)_\pm
$$
under the basic assumption $0<q<1<p$.

If $\lambda = 0$, equation $(P_\lambda)_\pm$ is called the Lane-Emden-Fowler equation and arises in the
boundary--layer theory of viscous fluids (see \cite{r33}).  This equation goes back to the paper by H.~Lane \cite{lane} in 1869 and is originally motivated by Lane's interest in computing both the temperature and the density of mass on the surface of the sun.
Equation $(P_\lambda)_\pm$ describes the behaviour of the density of a gas sphere in hydrostatic equilibrium and the index $p$, which is called the polytropic index in astrophysics, is related to the ratio of the specific heats of the gas.
Problem \plpm\ may be also viewed as a prototype of pattern formation in biology
and is related to the steady--state problem for a chemotactic aggregation model
introduced by Keller and Segel \cite{r1970}. Problem \plpm\  also plays an important role in
the study of activator-inhibitor systems modeling biological pattern formation, as
proposed by Gierer and Meinhardt \cite{r1972}. Problems of this type, as well as the associated
evolution equations, describe super-diffusivities phenomena. Such models have been proposed by de Gennes \cite{r12} to describe long range Van der Waals interactions in thin films spread on solid surfaces.
 This
equation also appears in the study of cellular automata and interacting particle systems
with self-organized criticality (see \cite{r5}), as well as to describe the flow over an impermeable
plate (see \cite{r3}). Problems of this type are obtained from evolution equations of the
form
$$u_{tt}=\mbox{div}\,(u^{m-1}\nabla u)+h(x,u)\qquad\mbox{in}\ \Omega\times (0,T)$$
through the implicit discretization in time arising in nonlinear semigroup theory (see \cite{r31}).

\medskip
Throughout this paper we assume that the variable weight functions $k,\, h\in L^\infty (\Omega)$ satisfy
$$\mbox{ess\,inf}_{x\in\Omega}k(x)>0\qquad\mbox{and}\qquad \mbox{ess\,inf}_{x\in\Omega}h(x)>0\,.$$

We are concerned with the existence of weak solutions of problems \plp \ and \plm, that is, functions $u\in H^1_0(\Omega)$ such that

(i) $u\geq 0$ a.e. on $\Omega$ and $u>0$ on a subset of $\Omega$ with positive measure;

(ii) for all $\varphi\in H^1_0(\Omega)$ the following identity holds
$$\intom \nabla u\nabla\varphi\,dx =\intom \left[\lambda k(x)u^q\pm h(x)u^p\right]\varphi\,dx\,.$$

The main results of this paper give a complete description of both cases arising in problem \plpm. The first theorem is concerned with problem \plp\ and establishes the
existence of a minimal solution, provided that $\lambda>0$ is small enough. This result gives a complete qualitative analysis of the problem and the proof combines monotonicity arguments with variational techniques. We refer to Bartsch and Willem \cite{bawi} for a related problem that is treated by means of a new critical point theorem, which guarantees the existence of infinitely many critical values of
an even functional in a bounded range. In this framework, Bartsch and Willem proved that there exists a sequence $(u_n)$ of solutions (not necessarily positive!) with corresponding energy ${\mathcal E}(u_n)<0$ and such that ${\mathcal E}(u_n)\ri 0$ as $n\ri\infty$.

\begin{thm}\label{th1}
Assume $0<q<1<p<2^\star -1$. Then there exists a positive number $\lambda^*$ such that the following properties hold:

a) for all $\lambda\in (0,\lambda^*)$ problem \plp \ has a minimal solution $u(\lambda)$. Moreover, the mapping $\lambda\longmapsto u(\lambda)$ is increasing.

b) problem \plp \ has a solution if $\lambda =\lambda^*$;

c) problem \plp\ does not have any solution if $\lambda>\lambda^*$.
\end{thm}

The next result is concerned with problem \plm\ and asserts that there is some $\lambda^*>0$ such that \plm\ has a nontrivial solution if $\lambda>\lambda^*$ and no solution exists provided that $\lambda <\lambda^*$.

\begin{thm}\label{th2}
Assume $0<q<1<p<2^\star -1$. Then there exists a positive number $\lambda^*$ the following properties hold:

a) if $\lambda>\lambda^*$, then problem \plm\ has at least one solution;



b) if $\lambda<\lambda^*$, then problem \plm\ does not have any solution.
\end{thm}

We point out that related results have been established by Il'yasov \cite{ilya} and Lubyshev \cite{lubi1,lubi2}. However, the proof techniques of \cite{ilya,lubi1,lubi2} are completely different from the arguments
in the present paper and they rely on the global fibering method introduced by S.I.~Pohozaev \cite{poho}, which is among the most powerful tools for studying nonlinear
differential equations. The fibering
scheme allows Il'yasov and Lubyshev to find constructively constrained minimization problems which
have the property of ground among all constrained minimization problems corresponding
to the given energy functional. We also refer to Bozhkov and Mitidieri \cite{bomi1,bomi2}
for relatively closed ``convex-concave" settings, both for nonlinear differential equations and for quasilinear elliptic systems with Dirichlet boundary condition.
 Garcia Azorero and Peral Alonso \cite{r112} used the mountain pass theorem to obtain the existence of a sign-changing solution in a related quasilinear setting.

\section{Proof of Theorem \ref{th1}}
We first prove that if $\lambda>0$ is sufficiently small then problem \plp\ has a solution. For this purpose we use the method of sub- and super-solutions.

Consider the problem
\begin{equation}\label{1sub}
\left\{\begin{array}{lll}
&\di -\Delta w=\lambda k(x)w^q, \qquad &\mbox{in}\ \Omega\\
&\di w>0, \qquad&\mbox{in}\ \Omega\\
&\di w=0, \qquad&\mbox{on}\ \partial\Omega\,.
\end{array}\right.
\end{equation}
Then, by \cite{brekam}, problem \eq{1sub} has a unique solution $w$. We prove that the function $\underline u:=\varepsilon w$ is a sub-solution of problem \plp\ provided that $\ep>0$ is small enough. For this purpose it suffices to show that
$$\ep\lambda k(x)w^q\leq\lambda k(x) \ep^qw^q+h(x)\ep^pw^p\qquad\mbox{in $\Omega$},$$
which is true for all $\ep\in (0,1)$.

Let $v$ be the unique solution of the linear problem
$$\left\{\begin{array}{lll}
&\di -\Delta v=1, \qquad &\mbox{in}\ \Omega\\
&\di v=0, \qquad&\mbox{on}\ \partial\Omega\,.
\end{array}\right.
$$
We prove that if $\lambda>0$ is small enough then there is $M>0$ such that $\overline u:=Mv$
is a super-solution of \plp. Therefore it suffices to show that
\begin{equation}\label{1super}
M\geq \lambda k(x)(Mv)^q+h(x)(Mu)^p\qquad\mbox{in $\Omega$}.
\end{equation}
Set
$$A:=\|k\|_{L^\infty}\|\cdot \|v\|_{L^\infty}^q\qquad\mbox{and}\qquad
B:=\|h\|_{L^\infty}\|\cdot \|v\|_{L^\infty}^p.$$
Thus, by \eq{1super}, it is enough to show that there is $M>0$ such that
$$M\geq\lambda AM^q+BM^p.$$
or, equivalently,
\begin{equation}\label{2super}1\geq\lambda AM^{q-1}+BM^{p-1}.\end{equation}
Consider the mapping $(0,\infty)\ni t\longmapsto \lambda At^{q-1}+Bt^{p-1}$. A straightforward computation shows that this function attains its minimum for $t=C\lambda^{1/(p-q)}$, where $C=[AB^{-1}(1-q)(p-1)^{-1}]^{(q-1)/(p-q)}$. Moreover, the global minimum of this mapping is $$(AC^{q-1}+BC^{p-1})\lambda^{(p-1)/(p-q)}.$$
This shows that condition \eq{2super} is fulfilled for all $\lambda\in (0,\lambda_0]$ and $M=
C\lambda^{1/(p-q)}$, where $\lambda_0>0$ satisfies
$$(AC^{q-1}+BC^{p-1})\lambda_0^{(p-1)/(p-q)}=1. $$

It remain to argue that $\ep w\leq Mv$. This is a consequence of the maximum principle (see \cite{pucser}), provided that $\ep>0$ is small enough. Thus, problem \plp\ has at least one solution $u(\lambda)$ for all $\lambda<\lambda^*$.

Set
$$\lambda^*:=\sup\{\lambda>0;\ \mbox{problem \plp\ has a solution}\}.$$
By the definition of $\lambda^*$, problem \plp\ does not have any solution if $\lambda>\lambda^*$.
In what follows we claim that $\lambda^*$ is finite. Denote
$$m:=\min\left\{\mbox{ess\,inf}_{x\in\Omega}k(x), \mbox{ess\,inf}_{x\in\Omega}h(x)\right\}>0.$$
Let $\lambda'>0$ be  such that $m(\lambda'+t^{p-q})>\lambda_1t^{1-q}$ for all $t\geq 0$, where $\lambda_1$ stands for the first eigenvalue of $(-\Delta)$ in $H^1_0(\Omega)$. Denote by $\varphi_1>0$ an eigenfunction of the Laplace operator corresponding to $\lambda_1$. Since $u(\lambda)$ solves \plp\ we have for all $\lambda >\lambda'$,
$$\begin{array}{ll}
\di \lambda_1\intom u(\lambda)\varphi_1dx&\di =\intom\left(\lambda k(x)u(\lambda)^q+h(x)u(\lambda)^p\right)\varphi_1dx\\
&\geq\di\intom m(\lambda u(\lambda)^q+u(\lambda)^p)\varphi_1dx>\lambda_1\intom u(\lambda)\varphi_1dx.\end{array}$$
This implies that $\lambda^* \leq\lambda'<+\infty$, which proves our claim.

Let us now prove that $u(\lambda)$ is a minimal solution of \plp. Consider the sequence $(u_n)_{n\geq 0}$ defined by $u_0=w$ ($w$ is the unique solution of \eq{1sub}) and $u_n$ is the unique solution of the problem
$$\left\{\begin{array}{lll}
&\di -\Delta u_n=\lambda k(x)u_{n-1}^{q}+h(x)u_{n-1}^{p}, \qquad &\mbox{in}\ \Omega\\
&\di u_n>0, \qquad&\mbox{in}\ \Omega\\
&\di u_n=0, \qquad&\mbox{on}\ \partial\Omega\,.
\end{array}\right.
$$
Then, by the maximum principle, $u_n\leq u_{n+1}\leq u(\lambda)$. Moreover, by the same argument as in \cite{ARC}, the sequence $(u_n)_{n\geq 0}$ converges to $u(\lambda)$. In order to show that $u(\lambda)$ is a minimal solution, let $U$ be an arbitrary solution of problem \plp. Thus, by the maximum principle, $w=u_0\leq U$ and, by recurrence, $u_n\leq U$ for all $n\geq 1$. It follows that $u(\lambda)\leq U$. At this stage it is easy to deduce that the mapping $\lambda\longmapsto u(\lambda)$ is increasing. Fix $0<\lambda_1<\lambda_2<\lambda^*$. Then $u(\lambda_2)$ is a super-solution of problem $(P_{\lambda_1})_+$, hence, by minimality, $u(\lambda_1)\leq u(\lambda_2)$. The fact that $\lambda_1<\lambda_2$, combined with the maximum principle implies that $u(\lambda_1)< u(\lambda_2)$.

It remains to show that problem \plp\ has a solution if $\lambda=\lambda^*$. For this purpose it is enough to prove that $(u(\lambda))$ is bounded in $\huo$ as $\lambda\ri\lambda^*$. Thus, up to a subsequence, $u(\lambda)\rightharpoonup u^*$ in $\huo$ as $\lambda\ri\lambda^*$, which implies that $u^*$ is a weak solution of \plp \ provided that $\lambda=\lambda^*$. Moreover, since the mapping $\lambda\longmapsto u(\lambda)$ is increasing, it follows that $u^*\geq 0$ a.e. on $\Omega$ and $u^*>0$ on a subset of $\Omega$ with positive measure. A key ingredient of the proof is that all solutions $u(\lambda)$ have negative energy. More precisely, if $\cale :\huo\ri\RR$ is defined by
$$\cale (u):=\frac12\intom |\nabla u|^2dx-\frac{\lambda}{q+1}\intom k(x)|u|^{q+1}dx-\frac{1}{p+1}\intom h(x)|u|^{p+1}dx$$
then \begin{equation}\label{ene}\cale (u(\lambda ))<0\qquad \mbox{for all $\lambda\in (0,\lambda^*)$}.\end{equation}
To deduce \eq{ene} we split the proof into the following steps:

(i) the solution $u(\lambda)$ is semi-stable, that is, the linearized operator $-\Delta -\lambda qk(x)u(\lambda)^{q-1}-ph(x)u(\lambda)^{p-1}$ is coercive:
$$\intom\left[|\nabla\psi|^2-(\lambda qk(x)u(\lambda)^{q-1}+ph(x)u(\lambda)^{p-1})\psi^2 \right]dx\geq 0\qquad\mbox{for all $\psi\in\huo$}.$$
Therefore
\begin{equation}\label{ene2}
\intom\left[|\nabla u(\lambda)|^2-(\lambda qk(x)u(\lambda)^{q+1}+ph(x)u(\lambda)^{p+1}) \right]dx\geq 0.
\end{equation}
This follows by the same proof as in \cite[Theorem 1.9]{hindawi}.

(ii) Since $u(\lambda)$ is a solution of \plp\ we have
\begin{equation}\label{ene3}\intom |\nabla u(\lambda)|^2dx=\lambda\intom k(x)  u(\lambda)^{q+1}dx+
\intom h(x)  u(\lambda)^{p+1}dx.\end{equation}

Combining relations \eq{ene2} and \eq{ene3} we deduce that
\begin{equation}\label{ene4}
\lambda (1-q)\intom k(x)  u(\lambda)^{q+1}dx\geq (p-1)\intom h(x)  u(\lambda)^{p+1}dx.\end{equation}

Next, we observe that relation \eq{ene3} implies
$$\begin{array}{ll}\di \cale (u(\lambda ))&\di =\lambda\left(\frac 12-\frac{1}{q+1}\right)\intom k(x)  u(\lambda)^{q+1}dx+
\left(\frac 12-\frac{1}{p+1}\right)\intom h(x)  u(\lambda)^{p+1}dx\\
&\di =-\lambda\,\frac{1-q}{2(q+1)}\intom k(x)  u(\lambda)^{q+1}dx+
\frac{p-1}{2(p+1)}\intom h(x)  u(\lambda)^{p+1}dx\\
&\di\leq -\lambda\,\frac{1-q}{2(q+1)}\intom k(x)  u(\lambda)^{q+1}dx+
\lambda\,\frac{1-q}{2(p+1)}\intom k(x)  u(\lambda)^{p+1}dx\leq 0,\end{array}$$
by \eq{ene4}.

To complete the proof, it remains to argue that $\sup_{\lambda<\lambda^*}\|u(\lambda)\|_{H^1_0}<+\infty$. This follows after combining relations \eq{ene} and \eq{ene2}, Sobolev embeddings, and using the fact that $k,\, h\in L^\infty(\Omega)$.
This completes the proof. \qed

\section{Proof of Theorem \ref{th2}}
The energy functional associated to problem \plm\ is $\calf :\huo\ri\RR$ and it is defined by
$$\calf (u):=\frac12\intom |\nabla u|^2dx-\frac{\lambda}{q+1}\intom k(x)|u|^{q+1}dx+\frac{1}{p+1}\intom h(x)|u|^{p+1}dx\,.$$
Set
$$\|u\|:=\left(\intom |\nabla u|^2dx\right)^{1/2};\
\|u\|_{q+1}:=\left(\intom | u|^{q+1}dx\right)^{1/(q+1)};\
\|u\|_{p+1}:=\left(\intom | u|^{p+1}dx\right)^{1/(p+1)}.$$

We first prove that $\calf$ is coercive. Indeed,
$$\calf (u)\geq\frac12\,\|u\|^2-C_1\|u\|_{q+1}^{q+1}+C_2\|u\|_{p+1}^{p+1}\,,$$
where $C_1=\lambda(q+1)^{-1}\|k\|_{L^\infty} $ and $C_2=(p+1)^{-1}\mbox{essinf}_{x\in\Omega}h(x)$ are positive constants. Since $q<p$, a straightforward computation shows that the mapping $(0,+\infty)\ni t\longmapsto At^{p+1}-Bt^{q+1}$ attains its global minimum $m<0$ at $$t=\left[\frac{B(q+1)}{A(p+1)}\right]^{1/(p-q)}\,.$$ Therefore
$$\calf (u)\geq \frac12\,\|u\|^2+m\,,$$
hence $\calf (u)\ri +\infty$ as $\|u\|\ri\infty$.

Let $(u_n)$ be a minimizing sequence of $\calf$ in $\huo$. Since $\calf$ is coercive, it follows that $(u_n)$ is bounded. Without loss of generality, we may assume that $u_n$ is non-negative and that $(u_n)$ converges weakly to some $u$ in $\huo$. Standard  arguments  based on the lower semi-continuity of the energy functional show that $u$ is a global minimizer of $\calf$, hence a non-negative solution of problem \plm.

In what follows we prove that the weak limit $u$ is a non-negative weak solution of problem \plm \ if $\lambda>0$ is large enough. We first observe that $\calf (0)=0$. So, in order to prove that the non-negative solution is nontrivial, it suffices to prove that there exists $\Lambda>0$ such that
$$\inf_{u\in\huo}\calf (u)<0\quad\mbox{for all $\lambda>\Lambda$}.$$
For this purpose we consider the constrained minimization problem
\begin{equation}\label{cop1}
\Lambda:=\inf\left\{\frac12\intom |\nabla v|^2dx+\frac{1}{p+1}\intom h(x)|v|^{p+1}dx;\ v\in\huo\ \mbox{and}\ \frac{1}{q+1}\intom k(x)|v|^{q+1}dx=1 \right\}.\end{equation}
Let $(v_n)$ be an arbitrary minimizing sequence for this problem. Then $(v_n)$ is bounded, hence we can assume that it  weakly converges to some $v\in\huo$ with
$$\frac{1}{q+1}\intom k(x)|v|^{q+1}dx=1\quad\mbox{and}\quad \Lambda:=\frac12\intom |\nabla v|^2dx+\frac{1}{p+1}\intom h(x)|v|^{p+1}dx\,.$$
Thus, $\calf (v)=\Lambda-\lambda<0$ for all $\lambda>\Lambda$.

Set
$$\lambda^*:=\inf\{\lambda>0;\ \mbox{problem \plm\ admits a nontrivial weak solution}\}\geq 0.$$

The above remarks show that $\Lambda\geq\lambda^*$ and that problem \plm\ has a solution for all $\lambda>\Lambda$. We now argue that problem \plm\  has a solution for {\it all} $\lambda>\lambda^*$. Fix $\lambda>\lambda^*$. By the definition of $\lambda^*$, there exists $\mu\in (\lambda^*,\lambda)$ such that
${\mathcal F}_\mu$ has a nontrivial critical point $u_\mu\in\huo$. Since $\mu<\lambda$, it follows that $u_\mu$ is a sub-solution of problem \plm. We now want to construct a super-solution that dominates $u_\mu$. For this purpose we consider the constrained minimization problem
\begin{equation}\label{cop2}
\inf\left\{\calf (v);\ v\in\huo\ \mbox{and}\ v\geq u_\mu \right\}.\end{equation}
The same arguments as those used to treat \eq{cop1} show that problem \eq{cop2} has a solution $u_\lambda\geq u_\mu$. Moreover, $u_\lambda$ is a solution of problem \plm, for all $\lambda>\lambda^*$. With the arguments developed in \cite[p.~712]{filpucrad} we deduce that problem \plm\ has a solution if $\lambda=\lambda^*$. The same monotonicity arguments as above show that \plm\  does not have any solution if $\lambda<\lambda^*$.

Fix $\lambda>\lambda^*$. It remains to argue that the non-negative weak solution $u$ is, in fact, positive.  Indeed, using similar arguments as in Pucci and Servadei \cite{pucser07}, which are based on the Moser iteration, we obtain that $u\in L^\infty(\Omega)$. Next, by bootstrap regularity, $u$ is a classical solution of problem \plm. Since $u$ is a non-negative smooth weak solution of the differential inequality $\Delta u -h(x)u^p\leq0$ in $\Omega$, with $p>1$, we deduce that $u>0$ in $\Omega$. This follows by applying the methods developed in Section 4.8 of Pucci and Serrin \cite{pucserhand} and the comments therein. This completes the proof. \qed

\medskip
{\bf Acknowledgments}. The authors are grateful the anonymous referees for the careful reading of the paper and their suggestions.
V.~R\u adulescu acknowledges the support through  Grant CNCSIS PCCE--8/2010
``Sisteme di\-fe\-ren\-\c{t}iale \^{\i}n analiza neliniar\u{a} \c{s}i aplica\c{t}ii".
D.~Repov\v{s} acknowledges the support by ARRS grant P1--0292--0101 ``Topology and geometry".

\end{document}